\documentclass[10pt]{article}

\usepackage{amsmath}
\usepackage{amsfonts}
\usepackage{amssymb}

\newtheorem{proposition}{Proposition}[section]
\newtheorem{theorem}{Theorem}[section]
\newtheorem{lemma}[theorem]{Lemma}
\newtheorem{corollary}[theorem]{Corollary}
\newtheorem{remark}[theorem]{Remark}

\newtheorem{definition}{Definition}

\newcommand{\braket}[2]{\langle #1,#2 \rangle}
\newcommand{\al}{\alpha}

\newcommand{\la}{\lambda}

\newcommand{\ra}{\rightarrow}
\def\phi{{\varphi}}

\DeclareSymbolFont{AMSb}{U}{msb}{m}{n}
\DeclareMathSymbol{\N}{\mathbin}{AMSb}{"4E}
\DeclareMathSymbol{\Z}{\mathbin}{AMSb}{"5A}
\DeclareMathSymbol{\R}{\mathbin}{AMSb}{"52}
\DeclareMathSymbol{\Q}{\mathbin}{AMSb}{"51}
\DeclareMathSymbol{\I}{\mathbin}{AMSb}{"49}
\DeclareMathSymbol{\C}{\mathbin}{AMSb}{"43}

\begin{document}

\title{Selfdual variational principles for periodic solutions of Hamiltonian and other dynamical systems}
\author{ Nassif  Ghoussoub\thanks{Partially supported by a grant
from the Natural Sciences and Engineering Research Council of Canada.  } \quad  and \quad Abbas Moameni\thanks{Research supported by a postdoctoral fellowship at the University of British Columbia.}
\\
\small Department of Mathematics,
\small University of British Columbia, \\
\small Vancouver BC Canada V6T 1Z2 \\
\small {\tt nassif@math.ubc.ca} \\
\small {\tt moameni@math.ubc.ca}
\\
}
\maketitle

\abstract
Selfdual variational principles are introduced in order to construct  solutions  for Hamiltonian and other dynamical systems which satisfy a variety of linear and nonlinear boundary conditions including many of the standard ones. These principles  lead to new variational  proofs of  the existence of parabolic flows with prescribed initial conditions, as well as periodic, anti-periodic and skew-periodic orbits of Hamiltonian systems. They are based on the theory of anti-selfdual Lagrangians introduced and developed recently  in \cite{G2}, \cite{G3} and \cite{G4}.

\section{Introduction}
The existence of a selfdual variational principle for gradient flows of convex functionals was conjectured in \cite{BE1} and established in \cite{GT1}. Similar selfdual variational principles were later introduced in \cite{GT3} and \cite{GM1}  for the resolution of certain gradient and Hamiltonian flows  that connect two prescribed Lagrangian submanifolds. In this paper, we introduce new anti-selfdual Lagrangians in order to construct variationally solutions of evolution equations that satisfy certain nonlinear boundary conditions. These include the more traditional ones, such as  the existence of flows with prescribed initial conditions, as well as periodic, anti-periodic and skew-periodic orbits. Our first variational principle typically deals with gradient flows of the form: 
\begin{equation}
-\dot x(t) = \partial\phi\big( t,x(t)\big)
\end{equation}
where $\phi (t, \,)$ is a convex lower semi-continuous function on a Hilbert space $H$. Our second   principle deals with Hamiltonian systems of the form:
\begin{equation}
-J\dot x (t) \in\partial\phi (t,x(t))
\end{equation}
where here $\phi (t, \cdot)$ is a convex lower semi-continuous functional on $H\times H$, and $J$ is the symplectic operator defined as $J(p,q)=(-q,p)$. In both cases, the prescribed conditions can be quite general but they include as particular cases the following more traditional ones:
\begin{itemize}
\item an initial value problem:  $x(0) = x_0$.
\item a periodic orbit:  $x(0) = x(T)$, 
\item an anti-periodic orbit: $ x(0) = -x(T)$ or
\item a skew-periodic orbit (in the case of a Hamiltonian system): $x(0)=Jx(T)$.
 \end{itemize}
 
 We are looking here for selfdual variational principles, and these  depend closely on the scalar product of the underlying path space. The novelty here is in the introduction of appropriate boundary Lagrangians $G$  which,  together with the main Lagrangian $L(t, x,p)$, yields an anti-selfdual Lagrangian on a path space equipped with an adequately  defined scalar product. The following space (scalar product) seems to be well adapted to our framework.
  
  Let  $[0,T]$ be a fixed real interval,  and let $L_H^2$ be the classical space of Bochner integrable functions from $[0,T]$  to  $H$. We consider the Hilbert space 
$A_H^2:=\left\{ u:[0,T]\ra H;\dot u\in L_H^2\right\}$
consisting of all absolutely continuous arcs $u:[0,T]\ra H$ equipped with the norm
\begin{eqnarray*}
\| u\|_{A_H^2}={\left\{ {\big\|\frac{u(0)+u(T)}{2}\big\|^2}_H
+\int_0^T \|\dot u\|_H^2 \, dt\right\} }^{\frac{1}{2}}
\end{eqnarray*}
We now recall the concept of anti-selfduality introduced in \cite{G2}. 
   \begin{definition} \rm Given a reflexive Banach space $X$, we say that a  convex lower semi-continuous function $L:X\times X^*\to \R \cup \{+\infty\}$ is an {\it anti-selfdual Lagrangian} if 
   \[
  L^*(p,x)=L(-x, -p) \quad \hbox{\rm for all $(x,p)\in X\times X^*$},
  \]
   where  here $L^*$ is the Legendre transform in both variables. 
   
 A  {\it  time dependent anti-selfdual Lagrangian} on $[0,T]\times X\times X^*$ is any function $L: [0,T]\times X\times X^*\to \R \cup \{+\infty\}$ that is  measurable with respect to   the  $\sigma$-field  generated by the products of Lebesgue sets in  $[0,T]$ and Borel sets in $X\times X^*$ and such that $L(t,\cdot, \cdot)$ is an anti-selfdual Lagrangian for every $t\in [0,T]$.
 
  The {\it Hamiltonian} $H_{L}$ of $L$ is the function defined
   on $[0,T]\times H\times H$ by:
   \[
   H_L(t, x,y)=\sup\{\langle y, p\rangle -L(t, x, p); p\in H\}
   \]
  \end{definition}
 Here is our first variational principle

\begin{theorem} \label{Principle.1} Consider a time dependent anti-selfdual Lagrangian 
$L(t,x,p)$ on $[0,T] \times H \times H$ where $H$ is a Hilbert space, and let $G$ be an anti-selfdual Lagrangian on $H\times H$. Consider on $A_H^2$ the following functional
\[
I(x)=\int_0^T L\big( t, x(t),\dot{x}(t)\big)\, dt
  +G\big(x(0)-x(T),
  \frac{x(0)+ x(T)}{2}\big). 
\]

 Assume the following conditions hold:
\begin{description}
\item[($A_1$)] \quad
$-\infty < \int_0^T L(t,x(t),0)\, dt\leq C\big(1+\|
x\|_{L^2_H}^2\big) ,\quad x\in L_H^2.$
\item[($A_2$)] \quad $\int_0^T H_L(t,0,x(t))\, dt\ra +\infty\quad\mbox{as}\quad \| x\|_{L^2_H} \ra +\infty.$
\item[($A_3$)] \quad $G$  is bounded from below and $0 \in {\rm Dom}_1 (G)$.
\end{description}

Then, there exists $\hat x\in A_H^2$ such that
\begin{eqnarray}
I(\hat x)&=&\inf\limits_{x\in A^2_H}I(x) =0 \label{P1}\\
 \big(-\dot{\hat{x}}(t),-\hat x(t)\big) &\in&\partial
   L\big(t,\hat x(t),\dot{\hat{x}}(t)\big)\quad {\rm for\,  all}\, t\in [0,T] \label{P2}\\
  \big(-\frac{\hat x(0)+\hat x(T)}{2}, \hat x (T)-\hat
x(0)\big)
 & \in&\partial G\big(\hat x(0)-\hat x(T),
  \frac{\hat x(0)+\hat x(T)}{2}\big). \label{P3}
\end{eqnarray}
\end{theorem}

The most basic time-dependent  anti-selfdual Lagrangians are of the
form
$
L(t,x,p)=\varphi (t, x) +\varphi^{*}(t, -p)
$
where for each $t$, the function $x\to \varphi (t, x)$ is convex and lower semi-continuous on $X$. Let  now $\psi: H \ra \R \cup\{+\infty\}$ be another convex lower semi-continuous function. The above principle then yields that  if  $-C\leq\int_0^T\phi\big( t,x(t)\big)\,dt\leq C\big(\|
x\|_{L_2^H}^2+1\big)$ and $\Phi (t,x):=\phi (t,x) +\frac{w}{2}|x|_H^2 +\langle f(t), x\rangle$, 
then the infimum of the functional 
\[
I(x)=\int_0^T \Phi( t, x(t))+\Phi^*(t, -\dot{x}(t))\, dt + \psi (x(0)-x(T)) +\psi^*(- \frac{x(0)+ x(T)}{2})
\]
on $A_H^2$ is zero and is attained at a solution $x(t)$ of the following equation
 \begin{eqnarray*}
 -\dot x(t) &=& \partial\phi\big( t,x(t)\big) +wx(t)+f(t) \quad {\rm for\,  all}\, t\in [0,T] \\
 -\frac{ x(0)+ x(T)}{2} & \in&\partial \psi (x(0)-x(T)).
\end{eqnarray*}
As to the various boundary conditions, we have to choose $\psi$ accordingly.  
  \begin{itemize}
\item Initial boundary condition $x(0)=x_0$ for a given $x_0\in
H$, then $\psi (x)= \frac{1}{4}\|x\|_H^2-\braket{x}{x_0}$.
\item Periodic solutions $x(0)=x(T)$, then $\psi$ is chosen as:
\begin{eqnarray*}\psi(x)=\left\{\begin{array}{ll}
0 \quad &x=0\\
+\infty &\mbox{elsewhere}.\end{array}\right.
\end{eqnarray*}
\item Anti periodic solutions $x(0)=-x(T)$, then $\psi (x)=0$ for each $x \in H.$
\end{itemize}
It is worth noting that while the main Lagrangian $L$ is expected to be smooth and hence its subdifferential coincides with its gradient --and the differential inclusion is often an equation, it is crucial that the boundary Lagrangian $G$ be allowed to be degenerate so as its subdifferential can cover the boundary conditions discussed above.
\bigskip

\noindent For the case of Hamiltonian systems we consider for simplicity $H=\R^N$ and let $X=H\times H$. We shall establish the following principle.
\begin{theorem} \label{Principle.2}
Let $\phi:[0,T]\times X\ra \R$ be  such that $(t,u)\ra\phi (t,u)$ is measurable in $t$ for each $u\in X$, and 
 convex and lower semi-continuous in $u$ for a.e. $t\in [0,T]$. Let $\psi: X\ra \R \cup \{\infty\}$   be convex and lower semi continuous on $X$ and assume the following conditions:
\begin{description}
\item ($B_1$)\quad There exists $\beta  \in (0,\frac{\pi}{2T})$ and  $\gamma, \al  \in L^2 (0,T;\mathbb{R_{+}})$ such that $- \al (t)\leq \phi(t,u)\leq\frac{\beta}{2} |u|^2+\gamma (t)$ for every
 $u\in H$ and all $t\in [0,T]$.
\item ($B_2$)\quad $\int_0^T \phi(t,u)\, dt\ra +\infty\quad\mbox{as}\quad | u| \ra +\infty.$
\item ($B_3$)\quad  $\psi$ is bounded from below and $0 \in Dom (\psi).$
\end{description}                         

(1) The infimum of the functional 
 \begin{eqnarray*}
 J_1(u)&=& \int_0^T \left[ \phi (t,u(t))+\phi ^*(t,-J\dot u (t))+\langle J\dot u(t),u(t)\rangle \right] dt\\
 &&\quad \quad \quad
 +\langle u(T)-u(0),J\frac{u(0)+u(T)}{2}\rangle + \psi \big( u(T)-u(0)\big) +\psi^*\big(
-J\frac{u(0)+u(T)}{2}\big)
\end{eqnarray*}
on $A_X^2$ is then equal to zero and is attained at a solution of 
\begin{eqnarray}\left\{ \begin{array}{lcl}
-J\dot u(t) &= &\partial \phi\big( t,u(t)\big),\\
-J\frac{u(T)+u(0)}{2} &\in &\partial \psi\big(
u(T)-u(0)\big).
\end{array}\right.
\end{eqnarray}
  
 (2)  The infimum of the functional 
\begin{eqnarray*}
J_2(u) =\int_0^T\left[\phi\big( t,u(t)\big)+\phi^*(t,-J\dot u (t))+ \langle J\dot u (t),u(t) \rangle \right]\, dt
   +\big( Ju(0),u(T)\big) +\psi\big( u(0)\big) +\psi^*\big( Ju(T)\big)
\end{eqnarray*}
on $A_X^2$ is also zero and is attained at a solution of 
\begin{eqnarray}\left\{ \begin{array}{lcl}
-J\dot u (t) &=&\partial\phi\big( t,u(t)\big),\\
Ju(T) &\in& \partial \psi\big( u(0)\big).
\end{array}\right.
\end{eqnarray}

 \end{theorem}
 In the applications, ${\psi}$ is to be chosen according to the required boundary conditions. For example: 
   \begin{itemize}
\item  Initial boundary condition $x(0)=x_0$ for a given $x_0\in H$. Use the functional $J_1$ with $  \bar {\phi} (t,x) =\phi (t, x-x_0)  $ and ${\psi}(x)=0$ at $0$ and $+\infty$ elsewhere.
 \item Periodic solutions $x(0)=x(T)$, or more generally $x(0)-x(T) \in K$ where $K$ is a closed convex subset of $H\times H$. Use the functional $J_1$ with ${\psi}$ chosen as:
\begin{eqnarray*} {\psi}(x)=\left\{\begin{array}{ll}
0 \quad &x\in K\\
+\infty &\mbox{elsewhere}.\end{array}\right.
\end{eqnarray*}
\item Anti-periodic solutions $x(0)=-x(T)$. Use the functional $J_1$ with ${\psi}(x)=0$ for each $x \in H.$
\item  Skew-periodic solutions $x(0)=Jx(T)$. Use the functional $J_2$ with ${\psi}(x)=\frac{1}{2}|x|^2$.
\end{itemize}
Section 2 deals with gradient flows and the proof of Theorem \ref{Principle.1}, while section 3 is concerned with Hamiltonian systems. This paper is self-contained but should be read in conjunction with \cite{G2}, \cite{G3} and \cite{GT1} which introduce selfduality and \cite{GM1} which deals with Hamiltonian systems that link Lagrangian submanifolds.
  
\section{Gradient flows with general boundary conditions}

\subsection{Anti-selfdual Lagrangians on path space}
We  now show how a boundary  anti-self dual Lagrangian allows us to ``lift" a time-dependent anti-selfdual Lagrangian to the path space $A_H^2$. Note that we can and will  identify the space $A_H^2$ with the product space $H\times L_H^2$, in such a way that its dual $(A_H^2)^*$ can also be identified with $H\times L_H^2$ via the formula
\begin{eqnarray*}
{\braket{u}{(p_1,p_0)}}_{A_H^2,H\times L_H^2}
  =\braket{\frac{u(0)+u(T)}{2}}{p_1} +\int_0^T
   \braket{\dot u(t)}{p_0(t)}\, dt
\end{eqnarray*}
where $u\in A_H^2$ and $(p_1,p_0(t))\in H\times L_H^2$.

\begin{proposition} \label{ASD}
Suppose $L$ is an anti-self dual Lagrangian on $[0,T]\times H\times H$ and that $G$ is an anti-selfdual Lagrangian on $H\times H$, then the Lagrangian defined on $A_H^2\times {(A_H^2)}^*=A_H^2\times (H\times L_H^2)$ by
\begin{eqnarray*}
{\cal M}(u,p)=\int_0^T L\big( t,u(t)+p_0(t),\dot u(t)\big)\, dt
  +G\big( u(0)-u(T)+p_1,\frac{u(0)+u(T)}{2}\big)
\end{eqnarray*}
is anti-self dual Lagrangian on $A_H^2\times (L_H^2\times H)$.
\end{proposition}

\paragraph{Proof:}
For $(q,v)\in A_H^2\times (A_H^2)^*$ with $q$ represented by $(q_0(t),q_1)$ we have
\begin{eqnarray*}
{\cal M}^*(q,v)&=&\sup\limits_{p_1\in H}\ \sup\limits_{p_0\in L_H^2}\ \sup\limits_{u\in A_H^2}
\Bigg\{\braket{p_1}{\frac{v(0)+v(T)}{2}}+\braket{q_1}{\frac{u(0)+u(T)}{2}}\\
& &\quad +\int_0^T \left[\braket{p_0(t)}{\dot v(t)}+\braket{q_0(t)}{\dot u}
  -L\big( t,u(t)+p_0(t),\dot u(t)\big)\right]\, dt\\
& &\quad -G\big( u(0)-u(T)+p_1,\frac{u(0)+u(T)}{2}\big)\Bigg\},
\end{eqnarray*}
making a substitution $u(0)-u(T)+p_1=a\in H$ and $u(t)+p_0(t)=y(t)\in L_H^2$ we obtain
\begin{eqnarray*}
{\cal M}^*(q,v) &=& \sup\limits_{a\in H}\ \sup\limits_{y\in L_H^2}\ \sup\limits_{u\in A_H^2}
\Bigg\{\braket{a+u(T)-u(0)}{\frac{v(0)+v(T)}{2}} +\braket{q_1}{\frac{u(0)+u(T)}{2}}\\
& &\quad +\int_0^T\left[\braket{y(t)-u(t)}{\dot v} +\braket{q_0(t)}{\dot u(t)}
  -L\big( t,y(t),\dot u(t)\big)\right]\, dt\\
& &\quad -G\big( a,\frac{u(0)+u(T)}{2}\big)\Bigg\}.
\end{eqnarray*}
Since $\dot u$ and $\dot v\in L_H^2$, we have:
$
\int_0^T\braket{u}{\dot v}=-\int_0^T\braket{\dot u}{v}+\braket{u(T)}{v(T)}
  -\braket{v(0)}{u(0)}
$
which implies
\begin{eqnarray*}
{\cal M}^*(q,v) &=& \sup\limits_{a\in H}\ \sup\limits_{y\in L_H^2}\ \sup\limits_{u\in A_H^2}
  \Bigg\{ \braket{a}{\frac{v(0)+v(T)}{2}}
  +\braket{u(T)}{\frac{v(0)+v(T)}{2}-v(T)}\\
& &\quad  +\braket{u(0)}{v(0)-\frac{v(0)+v(T)}{2}}
  +\braket{q_1}{\frac{u(0)+u(T)}{2}}\\
& &\quad +\int_0^T\left[\braket{y(t)}{\dot v}+\braket{\dot u(t)}{v(t)+q_0(t)}
  -L\big( t,y(t),\dot u(t)\big)\right]\, dt\\
& &\quad - G\big( a,\frac{u(0)+u(T)}{2}\big)\Bigg\} .
\end{eqnarray*}
Hence,
\begin{eqnarray*}
{\cal M}^*(q,v) &=& \sup\limits_{a\in H}\ \sup\limits_{y\in L_H^2}\ \sup\limits_{u\in A_H^2}
  \Bigg\{\braket{a}{\frac{v(0)+v(T)}{2}}
   +\braket{q_1+v(0)-v(T)}{\frac{u(0)+u(T)}{2}} -G\big( a,\frac{u(0)+u(T)}{2}\big)\\
& &\quad +\int_0^T\left[\braket{y(t)}{\dot v(t)}
  +\braket{\dot u(t)}{v(t)+q_0(t)}
  -L\big( t,y(t),\dot u(t)\big)\right]\, dt\Bigg\}.
\end{eqnarray*}
Identify now $A_H^2$ with $H\times L_H^2$ via the correspondence:
\begin{eqnarray*}
\big( b,f(t)\big) &\in & H\times L_H^2\longmapsto b+ \frac{1}{2}\left(
  \int_t^T f(s)\, ds-\int_0^t f(s)\, ds\right) \in A_H^2,\\
u &\in & A_H^2\longmapsto\big(\frac{u(0)+u(T)}{2},-\dot u(t)\big)
  \in H\times L_H^2.
\end{eqnarray*}
We finally obtain
\begin{eqnarray*}
{\cal M}^*(q,v) &=& \sup\limits_{a\in H}\ \sup\limits_{b\in H}
  \left\{\braket{a}{\frac{v(0)+v(T)}{2}}
   +\braket{q_1+v(0)-v(T)}{b} -G(a,b)\right\}\\
& &\quad +\sup\limits_{y\in L_H^2}\ \sup\limits_{r\in L_H^2}
  \left\{\int_0^T\left[\braket{y(t)}{\dot v(t)}
  +\braket{v(t)+q_0(t)}{r(t)}-L\big( t,y(t),r(t)\big)\right]\, dt\right\}\\
&=& G^*\big(\frac{v(0)+v(T)}{2},q_1+v(0)-v(T)\big)
   +\int_0^T L^*\big( t,\dot v(t),v(t)+q_0(t)\big)\, dt\\
&=& G\big( -q_1-v(0)+v(T),\frac{-v(0)-v(T)}{2}\big)
   +\int_0^T L\big( t,-v(t)-q_0(t),-\dot v(t)\big)\, dt\\
&=& {\cal M}(-v,-q).
\end{eqnarray*}

\subsection{Variational principles for gradient flows with general boundary conditions}

We now recall from \cite{G2} the following general result about minimizing anti-selfdual Lagrangians.

\begin{proposition}\label{nassif} Let ${\cal M}$ be a an anti-selfdual  Lagrangian on a reflexive Banach space $X\times X^{*}$ such  that for some $x_{0}\in X$, the function  $p\to {\cal M}(x_{0},p)$ is bounded above on  a neighborhood of the origin in $X^{*}$. Then there exists $\bar x\in X$, such that: 
 \begin{equation}
 \left\{ \begin{array}{lcl}
\label{eqn:existence}
  {\cal M}(\bar x, 0)&=&\inf\limits_{x\in X} {\cal M}(x,0)=0.\\
 \hfill (0, -\bar x) &\in & \partial {\cal M} (\bar x,0).
\end{array}\right.
 \end{equation}
 \end{proposition}

We can already deduce the following version of Theorem \ref{Principle.1} modulo a stronger hypothesis on the boundary Lagrangian.

 \begin{proposition} \label{Stringent.1}
Consider a time dependent anti-selfdual Lagrangian
$L(t,x,p) $ on $[0,T] \times H \times
H$ and an anti-selfdual lagrangian $G$ on $H\times H.$ Assume the
following conditions:
\begin{description}
\item ($A_1$) \quad $-\infty < \int_0^T L(t,x(t),0)\, dt\leq C\big(1+\|
x\|_{L^2_H}^2\big)$ for all $x\in L_H^2$. 
\item ($A_2$) \quad $G$  is bounded from below and
$
G(a,0)\leq C\big(\| a\|_H^2+1\big)$ for all $ a\in H$. 
 \end{description}
Then the functional
$
I(x)=\int_0^T L\big( t, x(t),\dot{x}(t)\big)\, dt
  +G\big(x(0)-x(T),
  \frac{x(0)+ x(T)}{2}\big) 
$
  attains its minimum at a path $\hat x\in A_H^2$ satisfying
\begin{eqnarray}
I(\hat x)&=&\inf\limits_{x\in A^2_H}I(x) =0\\
 \big(-\dot{\hat{x}}(t),-\hat x(t)\big) &\in&\partial
   L\big(t,\hat x(t),\dot{\hat{x}}(t)\big)\quad\forall t\in [0,T]\\
  \big(-\frac{\hat x(0)+\hat x(T)}{2}, \hat x (T)-\hat
x(0)\big)
 & \in&\partial G\big(\hat x(0)-\hat x(T),
  \frac{\hat x(0)+\hat x(T)}{2}\big) .
\end{eqnarray}

\end{proposition}

\paragraph{Proof:}
Apply Proposition \ref{nassif} to the  Lagrangian 
\begin{eqnarray*}
{\cal M}(u,p)=\int_0^T L\big( t,u(t)+p_0(t),\dot u(t)\big)\, dt
  +G\big( u(0)-u(T)+p_1,\frac{u(0)+u(T)}{2}\big)
\end{eqnarray*}
which is anti-selfdual on $A^2_H$ in view of Proposition \ref{ASD}. 
Noting that $I(x)={\cal M}(x,0)$, we  obtain  $\hat x(t)\in A_H^2$ such that
\begin{eqnarray*}
\int_0^T L\big( t,\hat x(t),\dot{\hat{x}}(t)\big)\, dt+G
  \left(\hat x(0)-\hat x(T),\frac{\hat x(0)+\hat x(T)}{2}\right) =0,
\end{eqnarray*}
which gives
\begin{eqnarray*}
0 &=&\int_0^T \left[L\big( t,\hat x(t),\dot{\hat{x}}(t)\big) + \braket
{\hat x(t)}{\dot{\hat{x}}(t)}\right]\, dt- \int_0^T \braket {\hat
x(t)}{\dot{\hat{x}}(t)} \,dt + G\big(\hat x(0)-\hat x(T),\frac{\hat x(0)+\hat x(T)}{2}\big) \\
&=&  \int_0^T \left[L\big( t,\hat x(t),\dot{\hat{x}}(t)\big) + \braket
{\hat x(t)}{\dot{\hat{x}}(t)}\right]\, dt- \frac{1}{2} | \hat  x(T) |^2+
\frac{1}{2} | \hat  x(0) |^2+ G
  \big(\hat x(0)-\hat x(T),\frac{\hat x(0)+\hat x(T)}{2}\big)\\
  &=&\int_0^T L\big( t,\hat x(t),\dot{\hat{x}}(t)\big) + \langle \hat
x(t),\dot{\hat{x}}(t)\rangle \, dt+ \braket {\hat x(0)-\hat
x(T)}{\frac{\hat x(0)+\hat x(T)}{2}}+ G
  \big(\hat x(0)-\hat x(T),\frac{\hat x(0)+\hat x(T)}{2}\big).
\end{eqnarray*}
Since $L(t,\cdot, \cdot)$ and $G$ are anti-selfdual Lagrangians we have
$
 L\big( t,\hat x(t),\dot{\hat{x}}(t)\big) + \braket
{\hat x(t)} {\dot{\hat{x}}(t)}\geq 0
$
and
\begin{eqnarray*}
 G
  \big(\hat x(0)-\hat x(T),\frac{\hat x(0)+\hat x(T)}{2}\big)+  \braket {\hat x(0)-\hat
x(T)}{\frac{\hat x(0)+\hat x(T)}{2}} \geq 0.
\end{eqnarray*}
 which means that 
$
 L\big( t,\hat x(t),\dot{\hat{x}}(t)\big) + \braket
{\hat x(t)} {\dot{\hat{x}}(t)}\,\big) = 0
$ for almost all $t\in [0,T]$, 
and
\begin{eqnarray*}
 G\big(\hat x(0)-\hat x(T),\frac{\hat x(0)+\hat x(T)}{2}\big)+  \braket {\hat x(0)-\hat
x(T)}{\frac{\hat x(0)+\hat x(T)}{2}} = 0.
\end{eqnarray*}
The result follows from the above identities and the limiting case in Fenchel-Legendre duality.$\square$\\

In order to complete the proof of Theorem \ref{Principle.1}, we need to perform an inf-convolution argument on the boundary Lagrangian $G$. We shall use the  following simple estimate
\begin{lemma}
Let $F:Y\mapsto
\mathbb{R}\cup\{\infty\}$ be a proper convex and lower semi
continues functional on a Banach space $Y$ such that $- \beta \leq F(y)\leq \frac{\al}{p}
\|y\|^{p}_{Y}+ \gamma$ with $\al >0, p >1, \beta\geq0,$ and $\gamma
\geq 0.$ Then for every $y^* \in \partial F(y)$ we have
\begin{eqnarray*}
\|y^*\|_{Y^*} \leq \Big \{ p \al^{\frac{q}{p}}( \|y\|_Y +\beta
+\gamma)+1 \Big \}^{p-1}.
\end{eqnarray*}
\end{lemma}
We shall also make frequent use of  the following  lemma \cite{G2}.

\begin{lemma} Let $G$ be an anti-selfdual Lagrangian on $X\times X^*$ and consider for each $\la>0$, its $\lambda$-regularization
\begin{eqnarray*}
G_\la (x,p):=\inf\left\{ G(z,p)+\frac{\| x-z\|^2}{2\la}+\frac{\la}{2}\| p\|^2; \, z\in X\right\} .
\end{eqnarray*}
Then,
\begin{enumerate}
\item $G_\la$ is also an anti-selfdual Lagrangian on $X\times X^*$ and
$
G_\la (x,0)\leq G(0,0)+\frac{\| x\|^2}{2\la}.
$
\item  If $(0,0)\in {\rm Dom}(G)$ and if 
 $x_\la\rightharpoonup x$ in $X$ and
$p_\la\rightharpoonup p$ weakly in $X^*$ and if $G
(x_\la, p_\la)$ is bounded from above, then
$
G(x,p) \leq \liminf\limits_{\la \ra 0} G_\la (x_\la, p_\la).
$
\end{enumerate}
\end{lemma}

\paragraph{Proof of Theorem \ref{Principle.1}:} Define for each $\la>0$, the Lagrangian $G_\la$ as in Lemma 2.2, 
 and apply Proposition \ref{Stringent.1} to obtain $x_\la \in A_H^2$ such that
\begin{eqnarray}
\int_0^T L\big( t, x_\la(t),\dot x_\la (t)\big)\, dt
  &+&G_\la \big( x_\la (0)-x_\la (T),\frac{x_\la (0)+x_\la (T)}{2}\big) =0\\
  \big(-\dot x_\la(t) ,-x_\la (t)\big)
  &\in&\partial L\big( t,x_\la (t),\dot x_\la (t)\big)\\
  \big(-\frac{x_\la (0)+x_\la (T)}{2},x_\la (T)-x_\la
(0)\big)
  &\in&\partial G_\la\big( x_\la (0)-x_\la (T),
  \frac{x_\la (0)+x_\la (T)}{2}\big).
\end{eqnarray}
We shall show that $ (x_\la)_\la$ is bounded in $A^2_H$.

For simplicity, we shall assume that $L$ has the form $L(t, x,p)=\phi (t, x) +\phi^*(t, -p)$. For such Lagrangians, Equation (13) yields  that 
$-\dot x_\la (t)=\partial_1 L\big( t,x_\la (t),0\big).
$
Multiply this equation by $x_\la (t)$ and integrate over $[0,T]\times\Omega$ to get
\begin{eqnarray*}
\int_0^T \braket {-\dot x_\la (t)}{ x_\la (t)} \, dt  =\int_0^T
  \braket {\partial_1 L (t,x_\la (t),0) }{x_\la (t)} \, dt,
\end{eqnarray*}
which  gives
\begin{equation}
-\frac{1}{2}| x_\la (T)|^2 +\frac{1}{2}| x_\la (0)|^2=\int_0^T
  \braket {\partial_1 L (t,x_\la (t),0) }{x_\la (t)} \, dt
   \geq\int_0^T H_L \big( t,0,x_\la (t)\big)\, dt.
\end{equation}
 Also, from (14) we have
 \begin{equation}
G_\la\big( x_\la (0)-x_\la (T),\frac{x_\la (0)+x_\la (T)}{2}\big)
 =-\braket { x_\la (0)-x_\la (T)}{
    \frac{x_\la (0)+x_\la (T)}{2}}\\
=\frac{1}{2}| x_\la (T)|^2 -\frac{1}{2}| x_\la (0)|^2.
\end{equation}

\noindent Combining~(15) and~(16) gives that
\begin{equation}
G_\la\big( x_\la (0)-x_\la (T),\frac{x_\la (0)+x_\la (T)}{2}\big)
  + \int_0^T H_L\big( t,0,x_\la (t)\big)\, dt
     \leq 0.
\end{equation} 
Since $G$ is bounded from below so is $G_\la$ which together with condition $(A_2)$ imply that
 $\int_0^T | x_\la (t)|^2\, dt$ is bounded.
 
 Now from condition $(A_1)$  and  the boundedness of $x_\la$ in $L_H^2$, we can apply 
Lemma 2.1 to get that  $- \dot x_\la (t)=\partial_1 L\big(
t,x_\la (t),0\big)$ is bounded in $L_H^2$. Hence, $x_\la$ is 
bounded in $A_H^2$, thus, up to a subsequence
$ x_\la (t)\rightharpoonup \hat x(t)$ in $A_H^2$,  $x_\la (0)\rightharpoonup \hat x(0)$ and $x_\la (T)\rightharpoonup  \hat x(T)$ in $H$. 

From (17), we have
$
G_\la\big( x_\la (0)-x_\la (T),\frac{x_\la (0)+x_\la (T)}{2}\big)\leq C, 
$
and  we obtain from Lemma 2.2 that 
\begin{eqnarray}
G\big( \hat x(0)-\hat x(T),\frac{\hat x(0)+\hat x(T)}{2}\big)
  \leq\liminf\limits_{\la \ra 0} G_\la
  \big( x_\la (0)-x_\la (T),\frac{x_\la (0)+x_\la (T)}{2}\big).
\end{eqnarray}
Now, if we let $\la\ra 0$ in (12),   then by considering (18) we get
\begin{eqnarray}
\int_0^T L\big( t,\hat x(t),\dot {\hat x}(t)\big)\, dt+
  G\big( \hat x(0)-\hat x(T),\frac{\hat x(0)+\hat x(T)}{2}\big)\leq 0.
\end{eqnarray}
On the other hand, for every $x \in A^2_H$ we have
\begin{eqnarray}
\int_0^T L\big( t,x(t),\dot x(t)\big)\, dt+
  G\big( x(0)-x(T),\frac{x(0)+x(T)}{2}\big)\ \geq 0
\end{eqnarray}
which  means $I(\hat x)=0$ and as in the
proof of Proposition \ref{Stringent.1}, $x(t)$  satisfies (\ref{P1}), (\ref{P2}),  and (\ref{P3}).
$\square$\\

The boundedness condition on $L$ may be too restrictive in applications, and one may want to replace 
the Hilbertian norm with a stronger Banach norm for which condition $(A1)$ is more likely to hold. For this situation,  we have the following result.
  
\begin{theorem}\label{Relax.1} Let $X\subset H \subset X^*$ be an evolution pair and let 
$\psi:[0,T] \times X \ra \mathbb{R} \cup \{+\infty\}$ be convex and
lower semi-continuous in $x \in X$ for a.e. $t \in [0,T]$ and
measurable in $t$ for every $x \in X.$ Consider the time-dependent
anti-selfdual Lagrangian, $L(t,x,p)=\psi (t,x)+\psi^* (t,-p)$ on
$[0,T] \times X \times X^*$ and an anti-selfdual Lagrangian $G$ on
$H\times H.$ Assume the following conditions:
\begin{description}
\item ($A'_1)$ \quad For some $p\geq 2$ and $C>0$, we have
 $-C\big(1+\| x\|_{L^p_X}^p\big) < \int_0^T L(t,x(t),0)\, dt\leq C\big(1+\| x\|_{L^p_X}^p\big)$ for every $x\in L^p_X$.
\item ($A'_2)$ \quad $G$ is bounded from below, $0 \in Dom (G)$ and for every $a\in H$, $G(a,b) \ra +\infty$ as $ \| b\|_{H}\ra +\infty. $
\end{description}

Then there exists $\hat x \in L^p_X$ with $\dot {\hat x} \in
L^q_{X^*}$ $(\frac{1}{p}+\frac{1}{q}=1)$, $\hat x(0), \hat x(T)
\in H$ and satisfying (\ref{P1}), (\ref{P2}),  and (\ref{P3}).

\end{theorem}
 {\bf Proof:} Here again we shall combine inf-convolution with Theorem \ref{Principle.1}.  For  $\la
>0$ consider the $\la-$regularization of $\psi$,   
\begin{equation}
\psi_\la (t,x)=\inf\limits_{y\in H}\left\{\psi (t,y)+\frac{|
x-y|_H^2}{2\la}\right\},
\end{equation}
where
\begin{eqnarray*}
\psi (t,y)=\left\{\begin{array}{ll}
\psi (t,y)\quad &y\in X\\
+\infty &y\in H-X.\end{array}\right.
\end{eqnarray*}
 Set  $L_\la(t, x,p)=\psi_\la (t,x)+\psi^*_\la (t,-p) .$
By Theorem \ref{Principle.1}, there exists $x_\la (t)\in A_H^2$ such that
\begin{eqnarray}
 \int_0^T L\big( t, x_\la(t),\dot x_\la (t)\big)\, dt
  &+&G \big( x_\la (0)-x_\la (T),\frac{x_\la (0)+x_\la (T)}{2}\big) =0\\
  \big(-\dot x_\la ,-x_\la (t)\big)
 & \in&\partial L\big( t,x_\la (t),\dot x_\la (t)\big)\\
 \big(-\frac{x_\la (0)+x_\la (T)}{2},x_\la (T)-x_\la
(0)\big)
  &\in&\partial G\big( x_\la (0)-x_\la (T),
  \frac{x_\la (0)+x_\la (T)}{2}\big)
\end{eqnarray}
We now  show that $(x_\la)_\la$ is bounded in an appropriate function space.  As in the proof of Theorem \ref{Principle.1}, we have
\begin{eqnarray}
G\big( x_\la (0)-x_\la (T),\frac{x_\la (0)+x_\la (T)}{2}\big)
  + \int_0^T H_{L_\la} \big( t, 0, x_\la (t)\big)\, dt
     \leq 0.
\end{eqnarray}
Since $\psi$ is convex and lower semi-continuous, there exists $i_\la
(x_\la)$ such that the infimum in (21) attains at $i_\la (x_\la),$
i.e.
\begin{eqnarray}
\psi_\la(t,x_\la)= \psi (t,i_\la (x_\la)) + \frac{\|x_\la-i_\la
(x_\la)\|^2}{2 \la}.
\end{eqnarray}
Therefore,
\begin{eqnarray}
 \int_0^T H_{L_\la} \big( t, 0, x_\la (t)\big)\, dt= \int_0^T H_L \big( t, 0, i_\la (x_\la (t))\big)dt+ \frac{\|x_\la-i_\la (x_\la)\|^2}{2 \la}\,
 dt.
\end{eqnarray}
Plug (27) in inequality (25) to get 
\begin{equation}
G\big( x_\la (0)-x_\la (T),\frac{x_\la (0)+x_\la (T)}{2}\big)
  + \int_0^T H_L \big( t, 0, i_\la (x_\la (t))\big) \, dt+ \frac{\|x_\la-i_\la (x_\la)\|^2}{2 \la}\, dt
     \leq 0.
 \end{equation}
By the coercivity assumptions in $(A'_1)$ , we obtain that $(i_\la (x_\la))_\la$
is bounded in $L^p(0,T; X)$ and $(x_\la)_\la$ is bounded in $ L^2(0,T; H)$.
It follows from  (23) and the structure of $L$ that
$
- \dot x_\la = \partial_1 L (t,i_\la (x_\la),0) ,
$
which together with the boundedness of $(i_\la (x_\la))_\la$ in  $L^p(0,T;
X)$, condition $(A'_1)$, and Lemma  2.1  imply that    $- (\dot x_\la)_\la$ is bounded
in $L^q(0,T; X^*)$. Also note that $x_{\la} (0)-x_{\la}(T)=
\int_{0}^{T}\dot {x_\la}(t)\, dt $  is therefore bounded in $X^*$. It follows from $(A'_2)$ that
$x_{\la} (0)+x_{\la}(T)$ is therefore bounded in $H$ and so is in $X^*$. Hence, up to a subsequence, we have
\begin{eqnarray}
i_\la (x_\la) \rightharpoonup \hat x \quad \text{  in  } L^p(0,T; X),\\
  \dot x_\la \rightharpoonup \dot
{ \hat x}\quad \text{  in  }  L^q(0,T; X^*),\\
{x_\la} \rightharpoonup  \hat x \quad \text{  in  }  L^2(0,T; H),\\
x_\la(0)\rightharpoonup \hat x(0) \quad \text{  in  }  X^*,\\
x_\la(T)\rightharpoonup \hat x(T) \quad \text{  in  } X^*.
\end{eqnarray}\\
On the other hand it follows from (22) and (26) that
\begin{equation}
G\big( x_\la (0)-x_\la (T),\frac{x_\la (0)+x_\la (T)}{2}\big)
  + \int_0^TL \big( t,i_\la (x_\la (t)),\dot x_\la \big)+ \frac{\|x_\la-i_\la (x_\la)\|^2}{2 \la} + \frac {\la}{2} \|\dot x_\la\|^2_H\, dt
    = 0.
\end{equation}
By letting $\lambda$ go to zero in (34), we get from (29)-(33) that
\[
 G\big( \hat x_\la (0)-\hat x_\la (T),\frac{\hat x_\la (0)+\hat x_\la (T)}{2}\big)
  + \int_0^TL \big( t, \hat x_\la (t)),\dot {\hat x}_\la \big) ^2_H\, dt
    \leq 0.
\]
It follows from $(A'_1)$ and the last inequality that $\hat x \in
L^p(0,T; X)$ and $\dot { \hat x} \in L^q(0,T; X^*).$ The rest of the
proof is similar to the proof of Proposition 2.3.

\begin{remark} \rm One can actually do without  the coercivity condition on $G$
in Theorem \ref{Relax.1}. Indeed, by using  the  $\la-$regularization $G_{\la}$
of $G$, we get the required coercivity condition on the second variable for
$G_{\la}$ and we obtain from Theorem \ref{Relax.1} that there exists  $
x_{\la} \in L^p(0,T; X)$ with  $\dot {x}_{\la} \in L^q(0,T; X^*)$
such that
\begin{eqnarray}
\int_0^T L\big( t, x_{\la}(t),\dot{x}_{\la}(t)\big)\, dt
  +G_{\la}\left( x_{\la}(0)- x_{\la}(T),
  \frac{ x_{\la}(0)+ x_{\la}(T)}{2}\right) = 0.
\end{eqnarray}
It follows from $(A_1)$ and the boundedness of $G_{\la}$ from below that $
(x_{\la})_\la$ is bounded in $ L^p(0,T; X)$, and since  $(\dot {x}_{\la})_\la$ is
bounded in $ L^q(0,T; X^*)$ this also means $(x_{\la}(0))_\la$ and $ (x_{\la}(T))_\la$ are
bounded in $H$. Hence, up to a subsequence we have
\begin{eqnarray}
 x_\la \rightharpoonup \hat x \quad \text{  in  } L^p(0,T; X),\\
  \dot x_\la \rightharpoonup \dot
{ \hat x}\quad \text{  in  }  L^q(0,T; X^*),\\
x_\la(0)\rightharpoonup \hat x(0) \quad \text{  in  } H,\\
x_\la(T)\rightharpoonup \hat x(T) \quad \text{  in  } H.
\end{eqnarray}
The rest of the proof is similar to the proof of Theorem \ref{Principle.1}. $\square$ 
\end{remark}

\subsection{Example}

As mentioned in the introduction, a typical example is 
\begin{eqnarray}
\left\{\begin{array}{rcl}
-\dot x(t) &= &\partial\phi\big( t,x(t)\big) +wx(t)+f(t)\\
x(0) &= &x_0\quad  \text{or} \quad x(0) = x(T) \quad \text{or} \quad
 x(0) = -x(T),\end{array}\right.
\end{eqnarray}
where $-C\leq\int_0^T\phi\big( t,x(t)\big)\,dt\leq C\big(\|
x\|_{L_2^H}^2+1\big)$ and $w >0$.

For the initial-value problem $x(0)=x_0$, we pick the boundary Lagrangian to be 
$G(x,p))= \frac{1}{4}|x|_H^2-\braket{x}{x_0}+|x_0-p|^2$, and so the associated functional becomes
 \[
I(x)=\int_0^T \Phi( t, x(t))+\Phi^*(t, -\dot{x}(t))\, dt + \frac{1}{4}|x(0)-x(T)|^2-\braket{x(0)-x(T)}{x_0}
+|x_0+ \frac{x(0)+ x(T)}{2}|^2 
\]
 where $\Phi (t,x):=\phi (t,x) +\frac{w}{2}|x|_H^2 +\langle f(t), x\rangle$, 
 The infimum of  $I$ on $A_H^2$ is zero and is attained at a solution $x(t)$ of the  equation. The boundary condition is then 
 \begin{eqnarray*}
-\frac{1}{2}(x(0)+x(T))=\partial_1 G\big(x(0)-
x(T), - \frac{x(0)+ x(T)}{2}\big)
   =\frac{1}{2}\big(x(0)- x(T)\big) -x_0,
\end{eqnarray*}
which gives that $x(0)=x_0.$

We can of course relax the conditions on $\phi$ by using again inf-convolution as was done in \cite{GT1} in the case where $\phi$ is autonomous, or as in Theorem 2.3.  

\section{Hamiltonian systems with general boundary conditions}

For a given Hilbert space $H$, we consider the subspace $H_{T}^1$  of   $A_{H}^2$ consisting of all periodic functions, equipped with the norm induced by $ A_H^2$. We also consider the  space $H_{-T}^1$ consisting of all
functions in $A_{H}^2$ which are anti-periodic, i.e.  $u(0)=-u(T).$
The norm of $H_{-T}^1$ is given by
$
\| u\|_{H_{-T}^1}=( \int_0^T |\dot u|^2 \, dt)^{\frac{1}{2}}.
$ We now  establish a few useful inequalities on $H_{-T}^1$, which can be seen as  the counterparts of  Wirtinger's inequality,
\begin{eqnarray*}
\int_0^T |u|^2 \, dt \leq \frac {T^2}{4 \pi^2} \int_0^T |\dot u|^2
\, dt  \quad {\rm for}\quad    u \in H_{T}^1   \text{   and  }\int_0^T
u(t) \, dt=0,
\end{eqnarray*}
and the Sobolev inequality on $H_{T}^1$,
\begin{eqnarray*}
 \|u\|^2_{\infty}  \leq \frac {T}{12} \int_0^T |\dot u|^2
\, dt  \quad {\rm for}\quad  u \in H_{T}^1   \text{   and  }\int_0^T
u(t) \, dt=0.
\end{eqnarray*}

\begin{proposition}
If $u \in  H_{-T}^1$ then
\begin{eqnarray}
\int_0^T |u|^2 \, dt \leq \frac {T^2}{ \pi^2} \int_0^T |\dot u|^2 \,
dt,
\end{eqnarray}
and
\begin{eqnarray}
 \|u\|^2_{\infty}  \leq \frac {T}{4} \int_0^T |\dot u|^2
\, dt.
\end{eqnarray}
\end{proposition}
\paragraph{Proof:} Since $u(0)=-u(T), $  $u$ has the Fourier
expansion of the form
$
u(t)= \sum_{k=-\infty}^{\infty}u_k \exp ((2k-1)i\pi t/T).
$
The Parseval equality implies that
\[
 \int_0^T |\dot u|^2 \,dt =\sum_{k=-\infty}^{\infty} T \big( (2k-1)^2
 \pi^2/T^2 \big) |u_k|^2  \geq
 \frac{\pi^2}{T^2}\sum_{k=-\infty}^{\infty}T |u_k|^2= \frac{\pi^2}{T^2}\int_0^T |u|^2
\, dt.
\]
The Cauchy-Schwarz inequality and the above  imply that for
$t \in [0,T],$
\begin{eqnarray*}
|u(t)|^2 &\leq &\left( \sum_{k=-\infty}^{\infty}|u_k| \right )^2 \\&
\leq & \left [  \sum_{k=-\infty}^{\infty} \frac{T}{\pi^2
(2k-1)^2}\right] \left [  \sum_{k=-\infty}^{\infty}  T \big( (2k-1)^2
\pi^2/T^2 \big)|u_k|^2 \right]\\ &=&\frac {T}{\pi^2}\sum_{k=-\infty}^{\infty}\frac {1}{(2k-1)^2}) \int_0^T |\dot u|^2
\, dt.
\end{eqnarray*}
and conclude by noting that 
 $\sum_{k=-\infty}^{\infty}\frac {1}{(2k-1)^2} =\frac{\pi^2}{4}$.
 
\begin{proposition}
Consider the space $A_X^2$ where $X= H \times H$ and let $J$ be the symplectic operator on $X$ defined as $J(p,q)=(-q,p)$.
\begin{enumerate}
\item If $H$ is any    Hilbert space, then  for every $u \in A_X^2$ 
\begin{eqnarray*}
\left|\int_0^T (J\dot u,u)\, dt+\left(J\frac{u(0)+u(T)}{2},u(T)-u(0)\right)\right| \leq\frac{T}{2}\int_0^T
\big|\dot u (t)\big|^2\, dt.
\end{eqnarray*}
\item  If $H$ is finite dimensional, then
\begin{eqnarray*}
\left|\int_0^T (J\dot u,u)\, dt+\left(
J\frac{u(0)+u(T)}{2},u(T)-u(0)\right)\right|
\leq\frac{T}{\pi}\int_0^T \big|\dot u (t)\big|^2\, dt.
\end{eqnarray*}
\end{enumerate}
\end{proposition}
\paragraph{Proof:} For part (i), note that each $u\in A_X^2$ can be written as follows,
\begin{eqnarray*}
u(t)=\frac{1}{2}\left( \int_0^t \dot u(s)\, ds-\int_t^T\dot u (s)\,
ds\right)+\frac{u(0)+u(T)}{2}.
\end{eqnarray*}
where $v(t)= u(t)-\frac {u(0)+u(T)}{2}=\frac{1}{2}\left( \int_0^t \dot u(s)\, ds-\int_t^T\dot u (s)\,
ds\right)$ clearly belongs to $H_{-T}^1$.
Multiplying both sides by $J\dot u$ and integrating over $[0,T]$, we get
\begin{eqnarray*}
\int_0^T \langle J\dot u,\, u\rangle\, dt=\frac{1}{2}\int_0^T\left\langle \int_0^t
\dot u(s)\, ds-\int_t^T \dot u(s)\, ds,\, J\dot u\right\rangle\,
dt+ \big\langle \frac{u(0)+u(T)}{2},\, \int_0^T J\dot u(t)\, dt\big\rangle
\end{eqnarray*}
Hence
\begin{eqnarray*}
\int_0^T \langle J\dot u,u\rangle\, dt-
\big\langle\frac{u(0)+u(T)}{2},\, J\big(u(T)-u(0)\big)\big\rangle=\frac{1}{2}\int_0^T\left\langle \int_0^t
\dot u(s)\, ds-\int_t^T \dot u(s)\, ds,\, J\dot u\right\rangle\,
dt\end{eqnarray*}
and since $J$ is skew-symmetric, we have
\begin{eqnarray}
\int_0^T \langle J\dot u,u\rangle \, dt+\big\langle
J\frac{u(0)+u(T)}{2},u(T)-u(0)\big\rangle
  =\frac{1}{2}\int_0^T\left\langle \int_0^t
\dot u(s)\, ds-\int_t^T \dot u(s)\, ds,\, J\dot u\right\rangle\,
dt
\end{eqnarray}
Applying H\"older's  inequality for the right hand side,
we get
\begin{eqnarray*}
\left|\int_0^T \langle J\dot u,u\rangle \, dt+\big\langle
J\frac{u(0)+u(T)}{2},u(T)-u(0)\big\rangle\right| \leq\frac{T}{2}\int_0^T
\big|\dot u (t)\big|^2\, dt.
\end{eqnarray*}
For part (ii), set $v(t)= u(t)-\frac {u(0)+u(T)}{2}$ and note that
\begin{eqnarray}
\int_0^T \langle J\dot u,u\rangle \, dt+\big\langle
J\frac{u(0)+u(T)}{2},u(T)-u(0)\big\rangle=\int_0^T (J\dot v,v)\, dt
\end{eqnarray}
Since $v \in H_{-T}^1$, H\"older's inequality  and  Proposition 3.1 imply,
\begin{eqnarray*}
\left |\int_0^T (J\dot v,v)\, dt  \right |& \leq& \left (\int_0^T
|v|^2 \, dt \right )^{\frac{1}{2}}\left (\int_0^T |J \dot v|^2 \, dt \right )^{\frac{1}{2}}\\
&\leq & \frac {T}{\pi}\left (\int_0^T |\dot v|^2 \, dt \right
)^{\frac{1}{2}}\left (\int_0^T |J \dot v|^2 \, dt \right
)^{\frac{1}{2}}\\ &=&\frac {T}{\pi} \int_0^T |\dot v|^2 \, dt
 =\frac {T}{\pi} \int_0^T |\dot u|^2 \, dt.
 \end{eqnarray*}
Combining this inequality with (44) yields  the claimed inequality.

\begin{proposition} If  $H= \mathbb{R}^N$ and $X= H \times H$, then 
the functional $F:A_X^2\ra \R$ defined by
\begin{eqnarray*}
F(u)=\int_0^T \langle J\dot u,u\rangle \, dt+\langle u(T)-u(0),J\frac{u(T)+u(0)}{2}\rangle
\end{eqnarray*}
is weakly continuous.
\end{proposition}

\paragraph{Proof:}
Let $u_k$ be a sequence in $A_X^2$ which converges weakly to $u$ in
$A_X^2$. The injection $A_X^2$ into $C([0,T];X)$ with natural norm $\|
\ \|_\infty$ is compact, hence $u_k\ra u$ strongly in $C([0,T];X)$ and
specifically $u_k(T)\ra u(T)$ and $u_k(0)\ra u(0)$ strongly in $X$.
Therefore
\begin{eqnarray}
\lim\limits_{k\ra +\infty} \left(
u_k(T)-u_k(0),J\frac{u_k(T)+u_k(0)}{2}\right)
  =\left( u(T)-u(0),J\frac{u(T)+u(0)}{2}\right)
\end{eqnarray}
Also, it is standard that $u\ra\int_0^T (J\dot u,u)\, dt$ is weakly
continuous (Proposition 1.2 in \cite{MW} ) which together with
(45) imply that $F$ is weakly continuous.

%
%
\subsection{A general variational principle for Hamiltonian systems}

In this section we establish Theorem \ref{Principle.2} under the assumption that  $H$ is finite dimensional  ($X=\R^{2N} $). We start with the following proposition which assumes a stronger condition on the boundary Lagrangian.
 
\begin{proposition} \label{Relax.2}
Let $\phi:[0,T]\times X\ra \R$, such that $(t,u)\ra\phi (t,u)$ is measurable in $t$ for each $u\in X$, and is 
 convex and lower semi-continuous in $u$ for a.e. $t\in [0,T]$. Let $\psi: X\ra \R \cup \{\infty\}$ be convex and lower semi continuous and assume the following conditions:
\begin{description}
\item ($B_1$)\quad There exists $\beta  \in (0,\frac{\pi}{2T})$ and  $\gamma, \al  \in L^2 (0,T;\mathbb{R_{+}})$ such that $- \al (t)\leq \phi(t,u)\leq\frac{\beta}{2} |u|^2+\gamma (t)$ for every $u\in X$ and a.e. $t\in [0,T]$. 
\item ($B'_2$)\quad There exist positive constants   $ {\al}_{1}, \beta_{1}, \gamma_{1} \in \mathbb{R}$ such that, for every $u\in X$  one has 
$- \al_1\leq \psi(u)\leq\frac{\beta_1}{2} |u|^2+ \gamma_1$.
\end{description}
   
\noindent (1) The infimum of the functional 
 \begin{eqnarray*}
 J_1(u)&=& \int_0^T \left[ \phi (t,u(t))+\phi ^*(t,-J\dot u (t))+\langle J\dot u(t),u(t)\rangle \right] dt\\
 &&\quad \quad \quad
 +\langle u(T)-u(0),J\frac{u(0)+u(T)}{2}\rangle + \psi \big( u(T)-u(0)\big) +\psi^*\big(
-J\frac{u(0)+u(T)}{2}\big)
\end{eqnarray*}
on $A_X^2$ is then equal to zero and is attained at a solution of 
\begin{eqnarray}\left\{ \begin{array}{lcl}
-J\dot u(t) &= &\partial \phi\big( t,u(t)\big)\\
-J\frac{u(T)+u(0)}{2} &= &\partial \psi\big(
u(T)-u(0)\big).
\end{array}\right.
\end{eqnarray}
  
\noindent (2)  The infimum of the functional 
\begin{eqnarray*}
J_2(u) =\int_0^T\left[\phi\big( t,u(t)\big)+\phi^*(t,-J\dot u (t))+\langle J\dot u (t),u(t)\rangle \right]\, dt
   +\big( Ju(0),u(T)\big) +\psi\big( u(0)\big) +\psi^*\big( Ju(T)\big)
\end{eqnarray*}
on $A_X^2$ is also equal to zero and is attained at a solution of 
\begin{eqnarray}\left\{ \begin{array}{lcl}
-J\dot u (t) &=&\partial\phi\big( t,u(t)\big)\\
Ju(T) &=& \partial \psi\big( u(0)\big).
\end{array}\right.
\end{eqnarray}

\end{proposition}
 
 The proof requires a few preliminary  lemmas, but first and
anticipating that the conjugate $\phi^*$  and
$\psi^*$ may not be finite everywhere,  we
start by  replacing $\phi$ and  $\psi$ with the
perturbations such as  $\phi_{\epsilon}(t,u)= \frac
{\epsilon}{2}\|u\|^2 + \phi(t,u)$ and $ \psi_{\epsilon}(u)= \frac {\epsilon}{2}\|u\|^2 + \psi(u)$. It is then  clear that
\begin{eqnarray}
\frac {1}{2(\beta+\epsilon)}|u|^2-\gamma (t) \leq  \phi^*_{\epsilon}(t,u)\leq  \frac {1}{2 \epsilon}|u|^2+\alpha (t) ,
\end{eqnarray}
and
\begin{eqnarray}
\frac {1}{2(\beta_1+\epsilon)}|u|^2-\gamma_1 \leq \psi^*_{\epsilon}(u)\leq  \frac {1}{2 \epsilon}|u|^2+\alpha_1 .
\end{eqnarray}
We now consider the Lagrangian ${\cal L}_{\epsilon}:A_X^2\times A_X^2\ra
\R $ defined  by
\begin{eqnarray*}
{\cal L}_{\epsilon}(v;u)&=&\int_0^T\left[ -\langle J\dot v (t) ,u(t)\rangle+\phi_{\epsilon} ^*(t,-J\dot u (t))-\phi_{\epsilon} ^*(t,-J\dot v
(t))+
\langle J\dot u(t),u(t)\rangle \right]\, dt\\
& &\quad +\left\langle u(T)-u(0),J\frac{u(T)+u(0)}{2}\right\rangle -\left\langle u(T)-u(0),J\frac{v(T)+v(0)}{2}\right\rangle
\\
& &\quad +\psi_{\epsilon}^*\big( -J\frac{u(T)+u(0)}{2}\big)
-\psi_{\epsilon}^*\big( -J\frac{v(T)+v(0)}{2}\big)
\end{eqnarray*}
and
\begin{eqnarray*}
J^{\epsilon}_1(u):&=& \int_0^T \left[ \phi_{\epsilon} (t,u(t))+\phi_{\epsilon} ^*(t,-J\dot u(t))+\big( J\dot u(t),u(t)\big)\right]\, dt\\
&& \quad +\big\langle u(T)-u(0),J\frac{u(0)+u(T)}{2}\big\rangle +\psi_{\epsilon} \big( u(T)-u(0)\big) + \psi_{\epsilon}^*\big(
-J\frac{u(0)+u(T)}{2}\big)
\end{eqnarray*}
To simplify the notation we use $C$ as a general positive constant.

\begin{lemma}
For   every $u\in A_X^2$, we have  $J_1(u)\geq 0$ and $J_1^\epsilon(u) \geq 0$. 
\end{lemma}

\paragraph{Proof:}
By the definition of   Legendre-Fenchel duality, one has
\begin{eqnarray*}
\phi (t,u(t))+\phi ^*(t,-J\dot u(t))+\langle J\dot
u(t),u(t)\rangle \geq 0 \quad  \quad \text{for }t \in [0,T],
\end{eqnarray*}
and
\begin{eqnarray*}
\psi \big( u(T)-u(0)\big) + \psi^*\big(
-J\frac{u(0)+u(T)}{2} \big)+ \langle
u(T)-u(0),J\frac{u(0)+u(T)}{2}\rangle\geq 0,
\end{eqnarray*}
which means $J_1(u)\geq 0$. The same applies to $J^\epsilon_1$. 

\begin{lemma}
For every $u\in A_X^2$, we have
$
I_{\epsilon}(u)= \sup\limits_{v\in A_X^2}{\cal L}_{\epsilon}(v,u).
$
\end{lemma}

\paragraph{Proof:}
First recall that one can identify  $A_X^2$ with $X\times L_X^2$ via the correspondence:
\begin{eqnarray*}
\big( x,f(t)\big) &\in & X\times L_X^2\longmapsto x+\frac{1}{2}\left(
  \int_0^t f(s)\, ds-\int_t^T f(s)\, ds\right) \in A_X^2\\
u &\in & A_X^2\longmapsto\left(\frac{u(0)+u(T)}{2},\dot u(t)\right)
  \in X\times L_X^2
\end{eqnarray*}
Thus, for every $u\in A_X^2$, we can write
\begin{eqnarray*}
\sup\limits_{v\in A_X^2}{\cal L}_{\epsilon}(v;u) &= & \sup\limits_{v\in X\times L_X^2}{\cal L}_{\epsilon}(v,u)\\
&=& \sup\limits_{f\in L^2(0,T;X)}\sup\limits_{x\in X}\Big\{\int_0^T \big[\langle -J f(t),u(t)\rangle +\phi_{\epsilon} ^*(t,-J\dot u(t))
 -\phi_{\epsilon} ^*(t,-J f(t))+(J\dot u(t),u(t))\big]\, dt\Big\}\\
&&\quad \quad  + \big \langle u(T)-u(0),J\frac{u(T)+u(0)}{2}\big \rangle -\big\langle u(T)-u(0),Jx\big\rangle+ \psi_{\epsilon}^*\big( -J\frac{u(T)+u(0)}{2}\big) -\psi_{\epsilon}^*(-Jx)\\
&=& \sup\limits_{f\in L^2(0,T;X)}\Big\{\int_0^T \big[\langle -J f(t),u(t)\rangle +\phi_{\epsilon} ^*(t,-J\dot u(t))
 -\phi_{\epsilon} ^*(t,-J f(t))+(J\dot u(t),u(t))\big]\, dt\Big\}\\
&&\quad \quad  + \sup\limits_{x\in X}\Big\{ \big \langle u(T)-u(0),J\frac{u(T)+u(0)}{2}\big \rangle -\big\langle u(T)-u(0),Jx\big\rangle\\
&&\quad \quad  +\psi_{\epsilon}^*\big( -J\frac{u(T)+u(0)}{2}\big) - \psi_{\epsilon}^*(-J x)\Big\}\\
&=& \int_0^T\left[ \phi_{\epsilon} (t,u(t))+\phi_{\epsilon} ^*(t,-J\dot u(t))+(J\dot u(t),u(t))\right]\, dt\\
&&\quad \quad  + \big\langle u(T)-u(0),J\frac{u(T)+u(0)}{2}\big\rangle +\psi_{\epsilon}\big( u(T)-u(0)\big)+\psi_{\epsilon}^*\big( -J\frac{u(T)+u(0)}{2}\big) \\
&=&J_1^{\epsilon}(u).
\end{eqnarray*}

\begin{lemma}
Under the assumptions $(B_1)$ and $(B'_2)$, we have for each
$0<\epsilon < \frac {1}{2}(\frac {\pi}{T}-2 \beta)$  the
following coercivity condition
\begin{equation}
{\cal L}_{\epsilon}(0,u)\ra +\infty \quad {\rm when} \quad \| u\|_{A_X^2}\ra +\infty .
\end{equation}
\end{lemma}

\paragraph{Proof:}
From (48) and (49) 
and since $\int_0^T \phi^*(t,0)\, dt$ and $\psi^*(0)$ are
finite,  we get
\begin{eqnarray*}
{\cal L}_{\epsilon}(0,u) &\geq& \frac{1}{2(\beta+\epsilon)}\int_0^T |\dot
u(t)|^2\, dt +\int_0^T \langle J\dot u(t),u(t)\rangle\, dt +\langle J\frac{u(0)+u(T)}{2}, u(T)-u(0)\rangle \\
& &\quad +\frac{1}{2(\beta_1+\epsilon)}{\left |\frac{u(0)+u(T)}{2}\right|}^2+C,
\end{eqnarray*}
 where $C$ is a constant. From part (ii) of Proposition 3.2, we have
\begin{eqnarray*}
\left|\int_0^T \langle J\dot u(t) ,u (t)\rangle \, dt+ \langle u(T)-u (0),
  J\frac{u(T)+u (0)}{2}\rangle \right|
  \leq\frac{T}{\pi}\int_0^T |\dot u (t)|^2\, dt
\end{eqnarray*}
Hence, modulo a constant, we obtain
\begin{eqnarray*}
{\cal L}_{\epsilon}(0,u) \geq
\left(\frac{1}{2(\beta+\epsilon)}-\frac{T}{\pi}\right)\int_0^T
{|\dot u(t)|}^2\, dt+\frac{1}{2(\beta_1+\epsilon)}
  {\left|\frac{u(0)+u(T)}{2}\right|}^2.
\end{eqnarray*}
Since $0<\epsilon < \frac {1}{2}(\frac {\pi}{T}-2 \beta)$, it follows 
that  $\frac{1}{2(\beta+\epsilon)}-\frac{T}{\pi}>0$ and  
${\cal L}_{\epsilon}(0,u)\ra +\infty$ as $\| u\|_{A_X^2}\ra +\infty$.  $\square $

\noindent Proposition 3.4 is now a consequence of the following Ky-Fan type min-max theorem which is essentially due to Brezis-Nirenberg-Stampachia (see \cite{BNS}).
 \begin{lemma} \label{KF} Let $Y$ be a  a reflexive Banach space  and let ${\cal L}(x,y)$ be a real valued function  on $Y\times Y  $ that satisfies the following conditions:
\begin{description}
\item (1) ${\cal L}(x,x) \leq 0$ for every $x\in Y$.
\item (2) For each $x\in Y$,  the function  $y \to {\cal L}(x,y)$ is concave.
\item (3) For each $y\in Y$, the function $x\to {\cal L}(x,y)$ is weakly lower semi-continuous.
\item (4) The set $Y_0=\{x\in Y; {\cal L}(x,0)\leq 0\}$ is bounded in $Y$.
\end{description}
Then there exists $x_{0}\in Y$ such that $\sup\limits_{y\in Y}{\cal L} (x_{0},y)\leq 0$.
\end{lemma}

\paragraph{Proof of Proposition \ref{Relax.2}:}
Let $0<\delta < \frac {1}{2}(\frac {\pi}{T}-2 \beta)$ and $0<
\epsilon < \delta.$ It is easy to see that the ${\cal L}_{\epsilon}:X\times
X\ra \R $ satisfies all the hypothesis of Lemma 3.4. It follows from
(48) and (49) that  ${\cal L}_{\epsilon}$ is finitely valued on $ X\times X $ and that 
for each $u\in X\times X$, ${\cal L}_{\epsilon}(u,u)=0$.
Lemma 3.3 gives that the set $Y=\{ u\in X, {\cal L}_{\epsilon}(0,u)\leq 0\}$ is
bounded in $X$. Moreover, for every $u\in X$,
 the function $v\ra {\cal L}_{\epsilon}(v,u)$ is concave and for every $v\in X$, $u\ra {\cal L}_{\epsilon}(u,v)$ is weakly
 lower semi-continuous by Proposition 3.3.
It follows that there exists $u_{\epsilon}\in X$ such that
$
I_{\epsilon}(u_{\epsilon})\leq\sup\limits_{v\in A_X^2}{\cal L}_{\epsilon}(v,u_{\epsilon})\leq 0.
$\\
In view of Lemma 3.1, we then have 
$I_{\epsilon}(u_{\epsilon})= 0$ which yields:
\begin{eqnarray}
I_{\epsilon}(u_{\epsilon}) &=& \int_0^T\left[\phi_{\epsilon}
(t,u_{\epsilon}(t))
+\phi_{\epsilon}^*(t,-J\dot{u_{\epsilon}}(t))+ \langle u_{\epsilon}(t),J\dot{u_{\epsilon}}(t)\rangle \right]\, dt \nonumber \\
&& \quad +\psi_{\epsilon} \big( u_{\epsilon}(T)-u_{\epsilon}(0)\big)
+\psi_{\epsilon}^*\big(
-J\frac{u_{\epsilon}(0)+u_{\epsilon}(T)}{2}\big)\nonumber\\
&&
  \quad + \langle u_{\epsilon}(T)-u_{\epsilon}(0), J\frac{u_{\epsilon}(0)+u_{\epsilon}(T)}{2} \rangle \nonumber \\
&=& 0.
\end{eqnarray}
We shall show that $u_{\epsilon}$ is bounded in $X$. From
Proposition 3.2,  we have
\begin{eqnarray*}
\left|\int_0^T (J\dot u_{\epsilon} (t),u_{\epsilon}(t) )\, dt+\langle
u_\epsilon (T)-u_\epsilon (0),
  J\frac{u_\epsilon (T)+u_\epsilon (0)}{2} \rangle \right|
  \leq\frac{T}{\pi}\int_0^T |\dot u_\epsilon (t)|^2\, dt
\end{eqnarray*}
which together with (51), yield
\begin{eqnarray*}
\int_0^T\left[\phi_{\epsilon} (t,u_{\epsilon}(t))
+\phi_{\epsilon}^*(t,-J\dot{u_{\epsilon}}(t))\right]\, dt-\frac{T}{\pi}\int_0^T |\dot u_{\epsilon} (t)|^2\, dt  
+ \psi_{\epsilon} \big( u_{\epsilon}(T)-u_{\epsilon}(0)\big)
+\psi_{\epsilon}^*\big(
-J\frac{u_{\epsilon}(0)+u_{\epsilon}(T)}{2}\big) \leq 0.
\end{eqnarray*}
This inequality together with the facts that $\phi_{\epsilon}$
and $\psi_{\epsilon}$ are bounded from below and $\phi^*_{\epsilon}$ and $\psi^*_{\epsilon}$ satisfy inequalities
(48) and (49) respectively, guarantee  the existence of
 a constant $C> 0$ independent of $\epsilon$ such that
 \begin{eqnarray*}
 \left(\frac{1}{2 (\beta + \delta)}-\frac{T}{\pi}\right)\int_0^T {|\dot u_{\epsilon}(t)|}^2\, dt+\frac{1}{2(\beta_1+ \delta)}
  {\left|\frac{u_{\epsilon}(0)+u_{\epsilon}(T)}{2}\right|}^2 &\leq
  &\\
 \left(\frac{1}{2(\beta+\epsilon)}-\frac{T}{\pi}\right)\int_0^T {|\dot u_{\epsilon}(t)|}^2\, dt+\frac{1}{2(\beta_1+ \epsilon)}
  {\left|\frac{u_{\epsilon}(0)+u_{\epsilon}(T)}{2}\right|}^2 &\leq & C,
\end{eqnarray*}
which means  $(u_{\epsilon})_\epsilon$ is bounded in $A^2_X$ and so, up to a
subsequence , there exists a $\bar u \in A^2_X$ such that
$u_{\epsilon}\rightharpoonup \bar u$ in $A^2_X.$ It is   easily seen
that
\begin{eqnarray*}
\int_{0}^{T} \phi^*_\epsilon(t, \dot {u}_{\epsilon}(t)) \,
dt:=\inf\limits_{v \in  L^2(0,T;X)}\int_{0}^{T}\left [ \phi^*(t,v(t))+\frac{| \dot {u}_{\epsilon}(t)-v(t) |^2}{2
\epsilon}\right]dt
\end{eqnarray*}
and since $ \phi^*$ is convex and lower semi continuous, there
exists $v_{\epsilon} \in  L^2(0,T;X)$  such that this infimum
attains  at $v_{\epsilon},$  i.e.
\begin{eqnarray*}
\int_{0}^{T} \phi^*_\epsilon (t,\dot {u}_{\epsilon}(t)) \,
dt=\int_{0}^{T}\left [ \phi^*(t, v_{\epsilon}(t))+\frac{| \dot
{u}_{\epsilon (t)}-v_{\epsilon}(t)|^2}{2 \epsilon}\right]dt.
\end{eqnarray*}
It follows from the above and the boundedness of $(u_{\epsilon})_\epsilon$ in $A^2_X,$ that 
there exists $C>0$ independent of $\epsilon$ such that
\begin{eqnarray*}
\int_{0}^{T} \phi^*_\epsilon (t, \dot {u}_{\epsilon}(t)) \,
dt=\int_{0}^{T}\left [ \phi^*(t,v_{\epsilon}(t))+\frac{| \dot
{u}_{\epsilon}(t)-v_{\epsilon}(t)|^2}{2 \epsilon}\right]dt<C.
\end{eqnarray*}
Since $\phi^*$ is bounded from below, we have
$
\int_{0}^{T}| \dot {u}_{\epsilon}(t)-v_{\epsilon}(t)|^2dt<C \epsilon
$
which means $ v_{\epsilon} \rightharpoonup \dot {\bar{u}}$ in
$L^2(0,T;X).$ Hence
\begin{eqnarray}
\int_{0}^{T}\phi^*(t,\dot {{\bar u}}(t))\, dt&\leq&
\liminf\limits_{\epsilon \rightarrow 0} \int_{0}^{T}\phi^*(t,v_{\epsilon}(t))dt\nonumber\\
& \leq&  \liminf\limits_{\epsilon
\rightarrow 0} \int_{0}^{T} \left [ \phi^*(t,v_{\epsilon}(t))+\frac{| \dot
{u_{\epsilon}}(t)-v_{\epsilon}(t)|^2}{2 \epsilon}\right]dt \nonumber\\ 
& =&\liminf\limits_{\epsilon \rightarrow 0}\int_{0}^{T}\phi^*_{\epsilon} (t, \dot
{u}_{\epsilon}(t))dt.
\end{eqnarray}
 Also,
\begin{eqnarray}
\int_{0}^{T}\phi(t, {{\bar u}}(t))\, dt&\leq&
\liminf\limits_{\epsilon \rightarrow 0} \int_{0}^{T}\phi (t,u_{\epsilon}(t))\, dt \nonumber\\
&\leq&  \liminf\limits_{\epsilon
\rightarrow
0} \int_{0}^{T} \left [ \phi(t,u_{\epsilon}(t))+\frac{ \epsilon}{2 }|u_{\epsilon}(t)|^2\right]dt \nonumber\\
& =&\liminf\limits_{\epsilon \rightarrow 0}\int_{0}^{T}\phi_{\epsilon} (t,
{u_{\epsilon}}(t))dt.
\end{eqnarray}
 It follows from (52) and (53) that,
\begin{eqnarray}
\int_0^T \left[\phi (t, \bar u(t)) +\phi^*(t,-J \dot{
\bar u}(t)) \right] \, dt \leq \liminf\limits_{\epsilon \rightarrow
0} \int_0^T\left[\phi_{\epsilon} (t,u_{\epsilon}(t)) +\phi_{\epsilon}^*(t,-J\dot{u_{\epsilon}}(t))\right] \, dt.
\end{eqnarray}
By the same argument  we arrive at,
\begin{eqnarray*}
\psi \big( \bar u(T)- \bar u(0)\big) + \psi^*\big(
-J\frac{\bar u(0)+\bar u(T)}{2}\big) \leq \liminf\limits_{\epsilon
\rightarrow 0} \Big \{ \psi_{\epsilon} \big(
u_{\epsilon}(T)-u_{\epsilon}(0)\big) + \psi_{\epsilon}^*\big(
-J\frac{u_{\epsilon}(0)+u_{\epsilon}(T)}{2}\big)\Big \}
\end{eqnarray*}
Also, from Proposition 3.3, we have
\begin{eqnarray}
\lim\limits_{\epsilon
\rightarrow 0}\int_0^T\langle u_{\epsilon}(t),J\dot{u_{\epsilon}}(t)\rangle\, dt   &+& \langle u_{\epsilon}(T)-u_{\epsilon}(0), J \frac{u_{\epsilon}(0)+u_{\epsilon}(T)}{2}\rangle\nonumber\\ 
&=&\int_0^T\langle \bar u(t),J\dot{\bar u}(t)\rangle dt   + \langle \bar
u(T)-\bar u(0), J \frac{\bar u(0)+\bar
u(T)}{2} \rangle.
\end{eqnarray}
 Combining  the above yields
\begin{eqnarray*}
I(\bar u) &=& \int_0^T\left[\phi (t,\bar u(t))
+\phi^*(t,-J\dot{\bar{u}}(t))+(\bar{u}(t),J\dot{\bar{u}}(t))\right]\, dt\\
&+& \psi \big(\bar u(T)-\bar u(0)\big) +\psi^*\big( -J\frac{\bar u(0)+\bar
u(T)}{2}\big)
  +\langle \bar u(T)-\bar u(0), J \frac{\bar u(0)+\bar u(T)}{2}\rangle\\
& \leq & \liminf\limits_{\epsilon \rightarrow 0}\Big\{
\int_0^T\left[\phi_{\epsilon} (t,u_{\epsilon}(t))
+\phi_{\epsilon}^*(t,-J\dot{u_{\epsilon}}(t))+(u_{\epsilon}(t),J\dot{u_{\epsilon}}(t))\right]\, dt\\
&+& \psi_{\epsilon} \big( u_{\epsilon}(T)-u_{\epsilon}(0)\big)
+\psi_{\epsilon}^*\big(
-J\frac{u_{\epsilon}(0)+u_{\epsilon}(T)}{2}\big)
+ \langle u_{\epsilon}(T)-u_{\epsilon}(0), J\frac{u_{\epsilon}(0)+u_{\epsilon}(T)}{2}\rangle\Big\} \\
&=& \liminf\limits_{\epsilon \rightarrow
0}I_{\epsilon}(u_{\epsilon})=0.
\end{eqnarray*}
On the other hand  Lemma 3.1 implies that  $I(\bar u)\geq 0,$ which
means the latter is zero, i.e.
\begin{eqnarray*}
I(\bar u) &=& \int_0^T\left[\phi (t,\bar u)
+\phi^*(t,-J\dot{\bar{u}})+(\bar{u},J\dot{\bar{u}})\right]\, dt\\
&+&\psi\big(\bar u(T)-\bar u(0)\big) +\psi^*\big( -J\frac{\bar u(0)+\bar u(T)}{2}\big)
  +\langle \bar u(T)-\bar u(0), J\frac{\bar u(0)+\bar u(T)}{2} \rangle=0.
\end{eqnarray*}
The result now  follows from the following identities and from the
limiting case in Legendre-Fenchel duality.
\begin{eqnarray*}
\phi (t,\bar u (t))+\phi^*(t,-J\dot{\bar{u}}(t))+(\bar
u(t),J\dot{\bar{u}}(t))=0
\end{eqnarray*}
\begin{eqnarray*}
\psi\big( \bar u(T)-\bar u(0)\big) +\psi^*\big(-J\frac{\bar u(0)+\bar u(T)}{2}\big)
  + \langle \bar u(T)-\bar u(0), J\frac{\bar u(0)+\bar u(T)}{2}\rangle  =0.
\end{eqnarray*}
$\square$\\

 We shall now use Proposition \ref{Relax.2}  to prove Theorem \ref{Principle.2}. For that 
 we shall $\la-$regularize the convex functional $\psi$, then 
use assumption $B_2$ of Theorem \ref{Principle.2} to derive  
uniform bounds and ensure convergence in $A^2_X$ when $\la$
approaches to $0$.  First recall that if 
$\psi_\la (x)=\inf\limits_{y\in X}\left\{ \psi (y)+\frac{\| x-y\|_X^2}{2\la}\right\} $ then its conjugate  $\psi^*_\la$ is equal to $\psi^* (x)+\frac{\la| x|^2}{2}$, which means that if $G(x,p)$ is the anti-selfdual Lagrangian $G(x,p)=\psi (x)+\psi^*(-p)$, then its $\la$-regularization is nothing but $G_\lambda (x,p)=\psi_\la (x)+\psi^*_\la (-p)$.

\paragraph{Proof of Theorem \ref{Principle.2}, Part (1):} The functional $\psi_\la$ satisfies
the condition $(B'_2)$ of Proposition \ref{Relax.2}, hence for each $\la >0$
there exists a $u_\la\in A_X^2$, such that
\begin{eqnarray}
I_\la (u_\la )&:=&\int_0^T\left[\phi(t,u_\la (t) )+\phi^*(t,-J\dot u_\la (t))+\langle J\dot u_\la (t) ,u_\la (t) \rangle \right]\, dt \nonumber \\
& &\quad +\langle u_\la (T)-u_\la (0),J\frac{u_\la (T)+u_\la
(0)}{2}\rangle
  +\psi_\la \big( u_\la (T)-u_\la (0)\big) +\psi_\la ^*\big(-J\frac{u_\la (T)+u_\la (0)}{2}\big) \nonumber \\
&=& 0.
\end{eqnarray}
We shall show $u_\la$ is bounded in $A_X^2$. From Proposition 3.2 
we obtain
\begin{eqnarray*}
\left |\int_0^T \langle J\dot u_\la ,u_\la \rangle\, dt+\langle u_\la (T)-u_\la
(0),
  J\frac{u_\la (T)+u_\la (0)}{2}\rangle \right|
  \leq\frac{T}{\pi}\int_0^T |\dot u_\la (t)|^2\, dt
\end{eqnarray*}
which together with (48) and (56) imply
\begin{eqnarray}
\psi_\la \big( u_\la (T)-u_\la (0)\big) +\psi_\la
^*\big(-J\frac{u_\la (T)+u_\la (0)}{2}\big)+\int_0^T \phi (t,u_\la (t))\, dt+\left(\frac{1}{2\beta}-\frac{T}{\pi}\right)
  \int_0^T |\dot u_\la (t)|^2 \, dt \leq 0.
\end{eqnarray}
Since $\psi$ is bounded from below so is $\psi_\la$.
Also, $0\in \mbox{ Dom}(\psi)$ which means $\psi^*$ and
consequently $\psi_\la ^*$ is bounded from below. Therefore
it follows from (57) that:
\begin{eqnarray}
\psi_\la \big( u_\la (T)-u_\la (0)\big) +\psi_\la
^*\big(-J\frac{u_\la (T)+u_\la (0)}{2}\big) \leq C,
\end{eqnarray}
and
\begin{eqnarray}
\int_0^T \phi(t,u_\la )\,
dt+\left(\frac{1}{2\al}-\frac{T}{2}\right)
  \int_0^T |\dot u_\la (t)|^2 \leq C,
\end{eqnarray}
where $C>0$ is a positive constant. 
It follows from the assumption  $(B_1), (B_2)$  and (59) that   $| u_\la
(t)|$ and  $\int_0^T |\dot u_\la |^2\, dt$ are  bounded.
Consequently $u_\la$ is bounded in $A_X^2$ and so, up to a subsequence,
$u_\la \rightharpoonup\bar{u}$ in $A_X^2$.

It follows from (58) and Lemma 2.2  that
\begin{eqnarray}
\psi\big(\bar u(T)-\bar u(0)\big) +\psi^* \big( -J\frac{\bar u(T)+\bar u(0)}{2}\big)
  \leq\liminf\limits_\la \psi_\la\big( u_\la (T)-u_\la (0)\big)
  +\psi_\la ^*\big( -J\frac{u_\la (T)+\bar u_\la(0)}{2}\big).
\end{eqnarray}
Also, from Proposition 3.3, we have
\begin{eqnarray}
\inf\limits_{\la
\rightarrow 0}\int_0^T\langle u_{\la}(t),J\dot{u_{\la}}(t)\rangle\, dt   &+& \langle u_{\la}(T)-u_{\la}(0), J \frac{u_{\la}(0)+u_{\la}(T)}{2}\rangle\nonumber\\
&=&\int_0^T\langle \bar u (t),J\dot{\bar u} (t) \rangle\, dt   +\langle \bar u(T)-\bar u(0), J \frac{\bar u(0)+\bar u(T)}{2}\rangle.
\end{eqnarray}
 Now, taking into account (60) and (61),  by letting $\la\ra 0$ in
(56) we obtain,
\begin{eqnarray*}
I(\bar u) &=& \int_0^T \left[\phi (t,\bar u (t))
  +\phi^*\big( t,-J\dot{\bar{u}} (t) \big)
  + \langle J\dot {\bar{u}}(t),\bar u (t) \rangle \right]\, dt\\
& &\quad +\langle \bar u (T)-\bar u (0),J\frac{\bar u (T)+\bar u
(0)}{2}\rangle
  +\psi\big(\bar u(T)-\bar u(0)\big) +\psi^*\big( -J\frac{\bar u (T)+\bar u (0)}{2}\big)\\
&\leq &\liminf\limits_{\la
\rightarrow 0} \Big \{\int_0^T\left[\phi (t,u_\la (t))+\phi^*(t,-J\dot u_\la (t))+\langle J\dot u_\la (t),u_\la (t)\rangle \right]\, dt\\
& &\quad +\langle u_\la (T)-u_\la (0),J\frac{u_\la (T)+u_\la
(0)}{2}\rangle
  +\psi_\la \big( u_\la (T)-u_\la (0)\big) +\psi_\la ^*\big(-J\frac{u_\la (T)+u_\la (0)}{2}\big)\Big\}\\
&=& \liminf\limits_{\la \rightarrow 0}I_\la (u_\la )=0.
\end{eqnarray*}
From Lemma 3.1, $I(\bar u)\geq 0$, which means the latter is zero.
The result follows from the following identities and from the
limiting case in Legendre-Fenchel duality
\begin{eqnarray*}
\phi (t,\bar u)+\phi^*(t,-J\dot{\bar{u}})+ \langle \bar
u,J\dot{\bar{u}}\rangle =0.
\end{eqnarray*}
\begin{eqnarray*}
\psi\big( \bar u(T)-\bar u(0)\big) +\psi^*\big( -J\frac{\bar u(0)+\bar u(T)}{2}\big)
  + \langle \bar u(T)-\bar u(0),  J \frac{\bar u(0)+\bar u(T)}{2} \rangle =0.
\end{eqnarray*}

\paragraph{Proof of Part (2):}
Note first that $\langle J\frac{u(0)+u(T)}{2},u(T)-u(0) \rangle= \langle 
Ju(0),u(T) \rangle. 
$
The corresponding Lagrangian $L_{\epsilon}:X\times
X\ra\R $ is defined as follows
\begin{eqnarray*}
L_{\epsilon}(v,u) &=& \int_0^T\left[ (-J\dot v (t),u (t))+\phi_{\epsilon}^* (t,-J\dot u (t))
  -\phi_{\epsilon}^*(t,-J\dot v (t))+\langle J\dot u (t),u (t) \rangle \right]\, dt\\
& &\quad +\langle Ju(0),u(T)\rangle -\langle Ju(0),v(T)\rangle
  +\psi_{\epsilon}^*\big( Ju(T)\big) -\psi_{\epsilon}^*\big(
  Jv(T)\big),
\end{eqnarray*}
The rest of the proof is quite similar to Part (1) and is left to the interested reader. 

\subsection{Applications} 

As mentioned in the introduction, one can choose the boundary  Lagrangian $\psi$ appropriately to 
solve Hamiltonian systems of the form 
\begin{eqnarray*}
\left\{\begin{array}{l}
-J\dot u (t) \in\partial\phi (t,u (t))\\
\hbox{$u(0)=u_0$, or $u(T)-u(0)\in K$, or $u(T)=-u(0)$ or $u(T)=Ju(0)$.}
\end{array}\right.
\end{eqnarray*}

 One can also use the method to solve  second order systems with convex potential and with prescribed nonlinear boundary conditions such as: 
  \begin{equation}\left\{ \begin{array}{lcl}
\hfill -\ddot{q} (t) &=& \partial \phi \big( t,q(t)\big)\\
-\frac{q(0)+q(T)}{2} &= &\partial\psi_1 \big( \dot q(T)- \dot q(0)\big),\\
\hfill \frac{\dot q(0)+ \dot q(T)}{2} &=&\partial\psi_2 \big(
q(T)-q(0)\big)
\end{array}\right.
\end{equation}
and
  \begin{equation}\left\{ \begin{array}{lcl}
\hfill \ddot{q} (t) &=& \partial \phi \big( t,q(t)\big)\\
-q(T) &= &\partial\psi_1 \big( \dot q(0)\big),\\
\hfill \dot q(T) &=&\partial\psi_2 \big( q(0)\big) 
\end{array}\right.
\end{equation}
 where $\psi_1$ and $\psi_2$ are convex and lower semi
continuous. One can deduce the following
\begin{corollary}
Let $\phi: [0,T]\times H\ra \R$ be  such that $(t,q)\ra \phi (t,q)$ is measurable in $t$ for each $q\in H$,
 convex and lower semi-continuous in $q$ for a.e. $t\in [0,T]$, and let $\psi_i: H \ra \R \cup \{\infty\}$, $i=1,2$
   be convex and lower semi continuous on $H$.
  Assume that the following conditions:
\begin{description}
\item[$A_1$:] There exists $\beta  \in (0,\frac{\pi}{2T})$ and  $\gamma, \al  \in L^2 (0,T;\mathbb{R_{+}})$ such that $- \al (t)\leq \mathcal \phi (t,q)\leq\frac{\beta^2}{2} |q|^2+\gamma (t)$ for every
 $q\in H$ and a.e. $t\in [0,T]$.
 \item[$A_2$:] $\int_0^T \mathcal \phi(t,q)\, dt\ra +\infty\quad\mbox{as}\quad | q|\ra +\infty.$
\item[$A_3$:] $\psi_1$ and $\psi_2$  are  bounded from below and $0 \in {\rm Dom} (\psi_i)$ for $i=1,2.$
\end{description}
Then equations (62) and (63) have at least one solution in $A_H^2.$
\end{corollary}
\paragraph{Proof:}
 Define $\Psi: H \times H\ra \R \cup \{\infty\}$ by  $\Psi (p,q):=\psi_1(p)+\psi_2(q)$ and  $\Phi:[0,T] \times H \times H\ra \R$ by $\Phi (t,u):=\frac
{\beta}{2} |p|^2+ \frac {1}{\beta}\phi \big( t,q(t)\big)$
where $u=(p,q).$ It is easily seen that  $\Phi$ is convex and
lower semi continuous in $u$ and that 
\begin{eqnarray*}
\hbox{$- \al (t)\leq \Phi(t,u)\leq\frac{\beta}{2} |u|^2+\frac {\gamma
(t) }{\beta}$ and $\int_0^T \Phi(t,u)\, dt\ra +\infty\quad\mbox{as}\quad | u| \ra
+\infty.$}
\end{eqnarray*}
 Also, from $A_3$, the function $\Psi$ is bounded from below
and  $0 \in {\rm Dom} (\Psi).$ By Theorem \ref{Principle.2}, the infimum of the
functional
\begin{eqnarray*}
I(u):&=& \int_0^T \left[ \Phi (t,u (t))+\Phi^*(t,-J\dot u (t))+\langle J\dot u(t),u(t)\rangle \right]\, dt\\
&& \quad + \langle u(T)-u(0),J\frac{u(0)+u(T)}{2}\rangle +\Psi\big(
u(T)-u(0)\big) +\Psi^*\big( -J\frac{u(0)+u(T)}{2}\big),
\end{eqnarray*}
on $A_X^2$ is zero and is attained at a solution of
\begin{eqnarray*}
\left\{\begin{array}{l}
-J\dot u (t)\in\partial\Phi (t,u(t)),\\
-J \frac {u(T)+u(0)}{2}=\partial\Psi (u(T)-u(0)).\end{array}\right.
\end{eqnarray*}
Now if we rewrite this problem for $u=(p,q),$ we get
\begin{eqnarray*}
- \dot p(t) &=& \frac {1}{\beta}\partial \phi \big( t,q(t)\big),\\
\dot q(t)&= &\beta p(t),\\
-\frac{q(T)+q(0)}{2} &= &\partial\psi\big( p(T)- p(0)\big),\\
\frac{p(T)+ p(0)}{2} &= &\partial\psi\big( q(T)-q(0)\big),
\end{eqnarray*}
and hence $q \in A_H^2$ is a solution of (61).
  
 As in the case of Hamiltonian systems, one can then solve variationally  the differential
equation $ -\ddot{q} (t) = \partial \phi \big( t,q(t)\big)$  with any one of the following boundary conditions:
\begin{description}
\item[(i)] Periodic: $ \dot q(T)= \dot q(0)$ and $ q(T)=  q(0).$
\item[(ii)] Antiperiodic: $ \dot q(T)= -\dot q(0)$ and $ q(T)=-  q(0).$
\item[(iii)] Initial value condition: $q(0)=q_0$ and $\dot
q(0)=q_1$ for given $q_0, q_1 \in H.$
\end{description}

\end{document}